\documentclass[a4paper,twoside,
 11pt]{article}
\usepackage{amssymb,amsfonts,amsmath}
\begin{document}
\newtheorem{defin}{~~~~Definition}
\newtheorem{prop}{~~~~Proposition}
\newtheorem{remark}{~~~~Remark}
\newtheorem{cor}{~~~~Corollary}
\newtheorem{theor}{~~~~Theorem}
\newtheorem{lemma}{~~~~Lemma}
\newtheorem{ass}{~~~~Assumption}
\newtheorem{con}{~~~~Conjecture}
\newtheorem{concl}{~~~~Conclusion}
\renewcommand{\theequation}{\thesection.\arabic{equation}}

\title{
Complete systems of invariants for rank 1 curves in 
Lagrange Grassmannians} 
\date{}
\author{Igor Zelenko\thanks{S.I.S.S.A., Via Beirut 2-4,
34014, Trieste, Italy; email: zelenko@sissa.it}} \maketitle 

\begin{abstract}
Curves in Lagrange Grassmannians naturally appear when one 
studies intrinsically "the Jacobi equations for extremals", 
associated with control systems and geometric structures. 
In this way one reduces the problem of construction of the 
curvature-type invariants for these objects to the much 
more concrete problem of finding of invariants of curves in 
Lagrange Grassmannians w.r.t. the action of the linear 
Symplectic group. In the present paper we develop a new 
approach to differential geometry of so-called rank 1 
curves in Lagrange Grassmannian, i.e., the curves with 
velocities being rank one linear mappings (under the 
standard identification of the tangent space to a point of 
the Lagrange Grassmannian with an appropriate space of 
linear mappings). The curves of this class are associated 
with "the Jacobi equations for extremals", corresponding to 
control systems with scalar control and to rank 2 vector 
distributions. In particular, we construct the tuple of $m$ 
principal invariants, where $m$ is equal to half of 
dimension of the ambient linear symplectic space, such that 
for a given tuple of arbitrary $m$ smooth functions there 
exists the unique, up to a symplectic transformation, rank 
1 curve having this tuple, as the tuple of the principal 
invariants. This approach extends and essentially 
simplifies the results of \cite{jac1}, where only the 
uniqueness part was proved and in rather cumbersome way. It 
is based on the construction of the new canonical moving 
frame with the most simple structural equation. 
%
%
%
%
%
\end{abstract}

\section{Statement of the problem and the results}
\setcounter{equation}{0} \indent 

Let $W$ be a $2m$-dimensional linear space provided with a 
symplectic form $\sigma$. Recall that an $m$-dimensional 
subspace $\Lambda$ of $W$ is called Lagrangian, if 
$\sigma|_\Lambda=0$. Lagrange Grassmannian $L(W)$ of $W$ 
 is the  set of all Lagrangian subspaces of $W$.
The linear Symplectic group acts naturally on $L(W)$. 
Invariants of curves in Lagrange Grassmannian w.r.t. this 
action are called symplectic 
The motivation to study differential geometry of curves in 
Lagrange Grassmannians comes from optimal control problems: 
it turns out that to any extremal of rather general control 
systems one can assign a special curve in some Lagrange 
Grassmannian, called the Jacobi curve (see \cite{agrachev}, 
\cite{agrgam1}, and Introduction to \cite{jac1} for the 
details). Symplectic invariants of Jacobi curves produce 
curvature-type differential invariants for these control 
systems. 

The natural differential-geometric problem is to construct 
a complete system of symplectic invariants for curves in 
Lagrange Grassmannians, i.e., some set of invariants such 
that there exists the unique, up to a symplectic 
transformation, curve in Lagrange Grassmannian with the 
prescribed invariants from this set. Some methods for 
construction and calculation of symplectic invariants of  
curves in Lagrange Grassmannians (including invariants of 
unparametrized curves) were given in \cite{jac1} and 
\cite{jac2}. Also the problem of finding a complete system 
of symplectic invariants for the special class of the 
so-called rank 1 curves in Lagrange Grassmannians (see 
Definition \ref{rankd} below) were partially solved there. 
In the present paper we solve this problem for the 
mentioned class of curves completely by developing another 
approach for the construction of symplectic invariants. 

Now we will briefly describe some main constructions of 
\cite{jac1} in order to specify in what sense our problem 
was partially solved and why it was not solved completely 
there. Note only that some results of the present paper 
(for example, Theorem \ref{maintheor}) do not depend on 
these constructions and in our opinion they are interesting 
by themselves. 
The key tool, used in \cite{jac1} for construction of 
symplectic invariants of curves in Lagrange Grassmannians, 
is an {\it infinitesimal cross-ratio} 
 of two tangent vectors 
$V_0,V_1$ at two distinct points $\Lambda_0$, $\Lambda_1$ 
in $L(W)$. In order to define it recall first that the 
tangent space $T_\Lambda L(W)$ to the Lagrangian 
Grassmannian $L(W)$ at the point $\Lambda$ can be naturally 
identified with the space ${\rm Quad}(\Lambda)$ of all 
quadratic forms on linear space $\Lambda\subset W$ or with 
the space ${\rm Symm} (\Lambda)$ of self-adjoint linear 
mappings from the space $\Lambda$ to the dual space 
$\Lambda^*$. Namely, take a curve $\Lambda(t)\in L(W)$ with 
$\Lambda(0)=\Lambda$. Given some vector $l\in\Lambda$, take 
a curve $l(\cdot)$ in $W$ such that $l(t)\in \Lambda(t)$ 
for all $t$ and $l(0)=l$. Define the quadratic form 
\begin{equation}
\label{quad} 
q_{\Lambda(\cdot)}(l)=\sigma(\frac{d}{dt}l(0),l).
\end{equation}
 Using the 
fact that the spaces $\Lambda(t)$ are Lagrangian, 
it is easy to see that the form $q_{\Lambda(\cdot)}(l)$ 
depends only on $\frac{d}{dt}\Lambda(0)$. One can consider 
also the self-adjoint linear mapping from $\Lambda$ to 
$\Lambda^*$, corresponding to this quadratic form. So, we 
have the mappings from $T_\Lambda L(W)$ to the spaces ${\rm 
Quad}(\Lambda)$ and ${\rm Symm}(\Lambda)$. 
 A simple counting of dimensions shows
that these mappings are bijections and they define the 
required identifications. Below we use these 
identifications 
without a special mentioning. Besides, given two Lagrangian 
subspaces $\Lambda_0$ and $\Lambda_1$, which are 
transversal, i.e. $\Lambda_0\cap\Lambda_1=0$, the map 
$v\mapsto \sigma(v, \cdot), \quad v\in \Lambda_1$, defines 
the canonical isomorphism from $\Lambda_1$ to 
$\Lambda_0^*$, which will be denoted by ${\mathcal 
B}_{\Lambda_0, \Lambda_1}$. 

Now we are ready to define the infinitesimal cross-ratio 
 of two tangent vectors 
$V_0,V_1$ at two distinct points $\Lambda_0$, $\Lambda_1$ 
in $L(W)$, where $\Lambda_0\cap\Lambda_1=0$. By above, 
$V_0$ is the self-adjoint linear mapping from $\Lambda_0$ 
to $\Lambda_0^*$, while 
$V_1$ is the self-adjoint linear mapping from $\Lambda_1$ 
to $\Lambda_1^*$. But $\Lambda_0^*\cong\Lambda_1$ via 
$\Bigl({\mathcal B}_{\Lambda_0, \Lambda_1}\Bigr)^{-1}$, 
while $\Lambda_1^*\cong\Lambda_0$ via $\Bigl({\mathcal 
B}_{\Lambda_1, \Lambda_0}\Bigr)^{-1}$. Therefore 
\begin{equation}
\label{crinfi} 
[V_0,V_1]_{\Lambda_0,\Lambda_1}\stackrel{def}{=}\Bigl({\mathcal 
B}_{\Lambda_1, \Lambda_0}\Bigr)^{-1}\circ V_1\circ 
\Bigl({\mathcal B}_{\Lambda_0, \Lambda_1}\Bigr)^{-1}\circ 
V_0 
\end{equation}
is well defined linear operator from the space $\Lambda_0$ 
to itself. This operator will be called an {\it 
infinitesimal cross-ratio} of a pair $(V_0,V_1)\in 
T_{\Lambda_0}L(W)\times T_{\Lambda_1}L(W)$. Actually, this 
notion is an infinitesimal version of the {\it cross-ratio} 
of four points in $L(W)$, which in turn is the 
generalization of the classical cross-ratio of four points 
on the projective line (see \cite{jac1} or \cite{Agr15} for 
the details). 

By constructions, the infinitesimal cross-ratio is the 
symplectic invariant of two tangent vectors 
at two distinct points of $L(W)$. Now let us show how to 
use this notion for the construction of symplectic 
invariants of a smooth curve $t\mapsto \Lambda(t)$ in 
$L(W)$. Suppose that the curve $\Lambda(\cdot)$ satisfies 
at some point $\tau$ the following condition: There exists 
a natural number $s$ such that for any representative 
$\Lambda^s_\tau(\cdot)$ of the $s$-jet of $\Lambda(\cdot)$ 
at $\tau$, there exists $t$ such that 
$\Lambda^s_\tau(t)\cap\Lambda(\tau)=0$. In this case the 
curve $\Lambda(\cdot)$ is called {\it ample at $\tau$}. The 
curve $\Lambda(\cdot)$ is called {\it ample}, if the last 
condition hold at any point $\tau$ of its segment of 
definition. To clarify this definition, let us give it in 
some coordinates: 
Let $W\cong\mathbb R^m\times\mathbb R^m$. 
The curve $t\mapsto\{(x,S_tx):x\in\mathbb R^n\}$ is ample 
at $\tau$ if and only if the function $t\mapsto 
det(S_t-S_\tau)$ has zero of {\it finite order} $k(\tau)$ 
at $\tau$. The number $k(\tau)$ is called {\it the weight} 
of an ample curve $\Lambda(\cdot)$ at $\tau$. Obviously, 
$k(\tau)$ is an integer valued upper semicontinuous 
function of $\tau$. Therefore it is locally constant on the 
open dense subset of the segment of definition of the 
curve. Note that any analytic monotone curve in $L(W)$, the 
image of which is not a point, is either ample or becomes 
ample in the Lagrange Grassmannian of another symplectic 
space, which obtained from $W$ after an appropriate 
symplectic factorization by the common subspace of all 
subspaces $\Lambda(t)$ (see Lemma 2.1 of \cite{jac1}). The 
term " monotone curve" means that its velocities are either 
nonnegative definite quadratic forms at any point or 
nonpositive definite quadratic forms at any point. Any 
curve $\Lambda(\cdot)$ in $L(W)$ such that $\dot 
\Lambda(t)$ is nondegenerated quadratic form on 
$\Lambda(t)$ (nondegenerated self-adjoint mapping from 
$\Lambda(t)$ to $\Lambda(t)^*$) has constant weight equal 
to $m$ ($=\frac{1}{2}\dim W$). We call such curves {\it 
regular}. Note that the set of all curves with constant 
weight in $L(W)$ is much wider than the set of all regular 
curves in $L(W)$. For example, any ample curve of rank 1 in 
$L(W)$ (see Definition \ref{rankd} below) at a generic 
point has the weight equal to $m^2$. 

The following proposition shows how to extract symplectic 
invariants from the infinitesimal cross-ratio: 

\begin{prop}{\rm(see \cite{jac1}, Lemma 4.2)}
\label{crprop} If the curve $\Lambda:I\mapsto L(W)$ has the 
constant finite weight $k$ on the segment of the definition 
$I$ , then the following asymptotic holds 
\begin{equation}
\label{crexp} {\rm trace} \Bigr[ \dot \Lambda(t_0)\,\mid 
\,\dot\Lambda(t_1)\Bigl]_{\Lambda(t_0),\Lambda(t_1)}= 
-\frac{k}{(t_0-t_1)^2}-g_{_{\Lambda}}(t_0,t_1), 
\end{equation}
 where $g_{_{\Lambda}}(t_0,t_1)$ 
  is a smooth function in the neighborhood of diagonal 
  $\{(t,t)| t\in I\}$. 
\end{prop}
Let us give the coordinate expression for the function 
$g_{_{\Lambda}}(t_0,t_1)$: If $W\cong\mathbb 
R^m\times\mathbb R^m$ 
, and $\Lambda(t)=\{(x, S_t x):x\in \mathbb R^m\}$ 
then 
\begin{equation}
\label{gdet} 
g_{_{\Lambda}}(t_0,t_1)=\cfrac{\partial^2}{\partial 
t_0\partial t_1} 
\ln\left(\frac{\det(S_{t_0}-S_{t_1})}{(t_0-t_1)^k}\right) 
\end{equation}
(the proof of the last formula follows from \cite{jac1}, 
see relations (4.9),(4.11), and Lemma 4.2 there).

 The 
function $g_{_{\Lambda}}$ is a "generating function" for 
the symplectic invariants in the following sense: suppose 
that it 
has the following 
expansion in the formal Taylor series at the point $(t,t)$ 
: 
\begin{equation}
\label{expg} g(t_0,t)\approx\sum_{i=0}^\infty 
\beta_i(t)(t_0-t)^i, 
\end{equation}
then all coefficients $\beta_i(t)$ are symplectic 
invariants of the curve $\Lambda(\cdot)$. In particular, 
the first appearing in (\ref{expg}) coefficient 
$\beta_0(t)$ ($=g_{_{\Lambda}}(t,t)$) produces the Ricci 
curvature, if one calculates it for Jacobi curves of 
Riemannian geodesics. 

The natural questions are whether the function 
$g_{_{\Lambda}}$ contains all information about the curve 
$\Lambda(\cdot)$ and what tuple of coefficients in 
expansion (\ref{expg}) of the function $g_{_{\Lambda}}$ 
constitutes a complete system of invariants of $\Lambda$? 
These questions were investigated in \cite {jac1}, section 
7, for so-called rank 1 curves.

\begin{defin}
\label{rankd} We say that a curve $\Lambda(\cdot)$ in 
$L(W)$ has rank $r$ at a point $t$, if its velocity 
$\dot\Lambda(t)$ is the linear self-adjoint mapping from 
$\Lambda(t)$ to $\Lambda(t)^*$ of rank $r$. A curve 
$\Lambda(\cdot)$ is called a rank $r$ curve in $L(W)$, if 
it has rank $r$ at any point $t$. 
\end{defin}

The motivation to study curves in $L(W)$ of rank less than 
$\frac{1}{2} \dim W$ at any point (which at first glance 
looks as rather degenerated case) comes from the fact that 
Jacobi curves associated with extremals of control systems 
with $r$-dimensional control space and $n$-dimensional 
state space (where $r<n$) are curves of rank not greater 
than $r$ in Lagrange Grassmannian of symplectic spaces of 
dimension equal usually to $2(n-1)$ or $2(n-2)$ (see 
Introduction to \cite{jac1} for the details). In 
particular, rank 1 curves in Lagrange Grassmannians appear 
as Jacobi curves associated with extremals of control 
systems with scalar control and with so-called abnormal 
extremals of rank 2 vector distributions (subbundles of the 
tangent bundle, see \cite{zel1} for the details). The fact 
that Jacobi curves are ample corresponds to some kind of 
controllability of the corresponding control system and to 
complete nonholonomicity (nonintegrability) of the 
corresponding distribution. 
 
 The main results in 
\cite{jac1}, concerning rank 1 curves in $L(W)$, is the 
following 

\begin{theor}
\label{uniqtheor} {\rm (see \cite{jac1}, Theorem 2)} The 
tuple $\{\beta_{2i}(t)\}_{i=0}^{m-1}$, where $\beta_j(t)$ 
as in (\ref{expg}), determines the curve $\Lambda(t)$, 
of rank $1$ and the constant finite weight 
 uniquely, up to a symplectic transformation. 
\end{theor}

For any $0\leq i\leq m-1$ the function $\beta_{2i}(t)$ will 
be called the {\it $(i+1)$-th principal curvatures} of the 
rank 1 curve $\Lambda(t)$. Theorem \ref{uniqtheor} actually 
states the uniqueness of rank 1 curve with the prescribed 
principal curvatures. {\it But the result about existence 
of rank 1 curve with the prescribed principal curvatures 
was missing}. Now I will try to explain the reason for it. 
Recall that the basis $(E_1,\ldots,E_m,F_1,\ldots F_m)$ of 
$W$ is called {\it Darboux}, if 
\begin{equation}
\label{darboux} \sigma(E_i,E_j)=\sigma(F_i,F_j)=0, \quad 
\sigma(F_i,E_j)=\delta_{ij},\quad \, 1\leq i,j\leq m, 
\end{equation} 
where $\delta_{ij}$ is the Kronecker symbol. To prove 
Theorem \ref{uniqtheor} we have constructed in \cite{jac1} 
a special canonical moving Darboux's frame for the given 
rank 1 curve in $L(W)$. To construct this frame we have 
used first the natural affine structure on the set 
$\Lambda(\tau)^\pitchfork$ of all Lagrangian subspaces 
transversal to $\Lambda(\tau)$, secondly the expansion of 
$t\mapsto \Lambda(t)$ (considered as the curve in this 
affine space $\Lambda(\tau)^\pitchfork$ with a singularity 
at $t=\tau$) into the Laurent series at $t=\tau$, and 
finally the fact that the free term of this Laurent series 
is the well-defined Lagrangian subspace transversal to 
$\Lambda(\tau)$, the {\it derivative subspace} (see 
Appendix below for the details). The first disadvantage of 
the canonical frame from \cite{jac1} is that the number of 
nontrivial entries in the matrix of its structural equation 
is much greater than the number of the functional 
parameters in our equivalence problem (which is equal to 
dimension of $m\times m$ symmetric matrices of rank 1, 
i.e., to $m$). So, this matrix does not give automatically 
a complete system of invariants: we need to choose some of 
the entries and to prove that all other nontrivial entries 
can be expressed by the chosen ones. Besides, the entries 
in the matrix of the considered structural equation are 
expressed in some nontrivial way by the principal 
curvatures. So, in order to prove Theorem \ref{uniqtheor} 
we had to analyze these expressions, which was rather 
nontrivial task. 
Another disadvantage is that the canonical frame from 
\cite{jac1} is not determined in the explicit way by the 
matrix of its structural equation: Even if some frame 
satisfies the structural equation with some prescribed 
functions substituted instead of the appropriate principal 
curvatures, it is not clear a priori whether this frame is 
canonical for the curve. Therefore it is not clear a priori 
whether the prescribed functions are exactly the 
corresponding principal curvatures of the curve. 
This is the reason why using this frame we did not succeed 
to prove the existence of the curve with the prescribed 
tuple of principal curvatures (in Remark \ref{diffproof} 
below we indicate the main technical difficulty that we met 
in this way). 

In the present paper  
we solve positively the problem of existence of the rank 1 
curves of the constant weight with the prescribed tuple of 
principal curvatures. For this we introduce a new very 
natural canonical moving Darboux's frame for a rank 1 curve 
in $L(W)$, which is uniquely defined by the matrix of its 
structural equation. It allows to obtain a new tuple of 
principal curvatures for which the uniqueness an existence 
results follow automatically. Namely, we have the following 

\begin{theor}
\label{maintheor} Let $W$ be a $2m$-dimensional linear 
symplectic space. For the given rank $1$ curve 
$\Lambda(\cdot)$ of the constant finite weight in $L(W)$ 
there exists the moving Darboux frame 
$$\bigl(E_1(t),\ldots, E_m(t),F_1(t),\ldots, F_m(t)\bigr)$$ 
such that the following two conditions hold: 
\begin{enumerate}
\item $\Lambda(t)={\rm span}\bigl(E_1(t),\ldots, E_m(t)\bigr)$;


\item The moving frame $\bigl(E_1(t),\ldots, E_m(t),F_1(t),\ldots, F_m(t)\bigr)$
satisfies the following structural equation:
\begin{equation}
\label{structeq} \left\{\begin{aligned} 
~&E_i^\prime(t)=E_{i+1}(t),\quad 1\leq i\leq m-1\\ 
~&E_m^\prime(t)=\pm F_m(t)\\ 
~&F_1^\prime(t)=\lambda_m(t)E_1(t)\\ 
~&F_i^\prime(t)=\lambda_{m-i+1}(t)E_i(t)-F_{i-1}(t),\quad 
2\leq i\leq m 
\end{aligned}\right.
\end{equation}
(in the second equation of (\ref{structeq}) the sign $"+"$ 
appears if the quadratic form $\dot\Lambda(t)$ is 
nonnegative definite, while the sign $"-"$ appears if the 
quadratic form $\dot\Lambda(t)$ is nonpositive definite). 
\end{enumerate}
In addition, if the moving frame $\bigl(E_1(t),\ldots, 
E_m(t),F_1(t),\ldots, F_m(t)\bigr)$ satisfies conditions 
1-3, then the only frame, which is different from it and 
satisfies the same conditions, is $\bigl(-E_1(t),\ldots,\\
-E_m(t),-F_1(t),\ldots, -F_m(t)\bigr)$. 
\end{theor}

The proof of Theorem \ref{maintheor} will be given in 
section 2. It consists basically of three steps: 
First the condition of rank 1 and the constant weight 
allows to construct the curve of complete flags in $W$, 
associated with our curve (see Lemmas \ref{constrank}, 
\ref{cwlem}, and Remark \ref{flag} below), secondly the 
presence of the symplectic structure allows to normalize 
the one-dimensional subspaces of the flags, which in turn 
gives the canonical basis on each subspace $\Lambda(t)$ 
(see Lemma \ref{normlem} and formula (\ref{Ecan})), and 
finally we complete this basis to the moving Darboux frame 
with the "most simple" structural equation (having the 
maximal possible number of zero entries in the matrix 
corresponding to it). 

From the uniqueness of the frame in Theorem \ref{maintheor} 
it follows immediately that each function $\lambda_i(t)$ in 
its structural equation is a symplectic invariant. It will 
be called the {\it $i$th modified principal curvature} of 
the curve $\Lambda(\cdot)$. Also, as a direct consequence 
of Theorem \ref{maintheor} we obtain the following 
\begin{theor}
\label{uniex} For the given tuple of $m$ smooth functions 
$\{\rho_i\}_{i=1}^m$ there exists the unique rank 1 curve 
in the Lagrange Grassmannian $L(W)$ such that its $i$-th 
\underline{modified} principal curvature coincides with 
$\rho_i(t)$ for any $1\leq i\leq m$. \end{theor} In other 
words, the tuple of the modified principal curvatures, 
defined by the structural equation (\ref{structeq}), 
constitutes the complete system of symplectic invariants 
for rank 1 curves of the constant rank. Besides, we have 
\begin{prop}
\label{equivpr}
 The following relations between the tuple of the principal 
curvatures $\{\beta_{2i}(t)\}_{i=0}^{m-1}$ from 
(\ref{expg}) and the tuple of the modified principal 
curvatures $\{\lambda_i(t)\}_{i=1}^m$ from (\ref{structeq}) 
hold: 
\begin{equation}
\label{equivbl} 
\lambda_{i}(t)=C_i \beta_{2i-2} +\Phi_i(t),
\quad 1\leq i\leq m 
\end{equation} 
where $C_i$ are nonzero constants, any $\Phi_i(t)$ is some 
polynomial expression (over $\mathbb R$) without free term 
w.r.t. the functions $\beta_{2j}(t)$, $0\leq j\leq i-2$ and 
their derivatives. \end{prop} 

\begin{remark}
\label{mpoint} {\rm Actually, the fact that the constants 
$C_{i}$, appearing in (\ref{equivbl}), are nonzero follows 
from Theorem \ref{uniex}: Assuming the converse, take the 
smallest $\bar i$ such that $C_{\bar i}=0$. Then from all 
relations (\ref{equivbl}) with $1\leq i<\bar i$ it follows 
that $\lambda_{\bar i}(t)$ is some polynomial expression 
w.r.t. the functions $\lambda_i(t)$, $1\leq i<\bar i$ and 
their derivatives. But this contradicts the fact that 
 the functions 
$\lambda_i(t)$, $1\leq i\leq m$, are independent according 
to Theorem \ref{uniex}.} $\Box$ 
\end{remark} 

 The proof of Proposition \ref{equivpr} will 
be given in section 3. As a direct consequence of Theorem 
\ref{uniex} and the previous proposition we obtain the 
following extension of Theorem \ref{uniqtheor}: 

\begin{theor}
\label{uniex1} For the given tuple of $m$ smooth functions 
$\{\rho_i\}_{i=1}^m$ there exists the unique rank 1 curve 
in the Lagrange Grassmannian $L(W)$ such that its $i$-th 
principal curvature coincides with $\rho_i(t)$ for any 
$1\leq i\leq m$. \end{theor} In other words, the tuple of 
the principal curvatures, defined by expansion 
(\ref{expg}), constitutes the complete system of symplectic 
invariants of rank 1 curves of the constant rank, the proof 
of which was the original goal of the present paper. 

Note that a kind of the complete system of invariants for 
regular curves (i.e., with nondegenerated velocities 
$\dot\Lambda(t)$) was constructed in \cite{Agr15}. In the 
forthcoming paper we will use the ideology of the proof of 
Theorem \ref{maintheor} in order to construct a kind of 
complete system of symplectic invariants for generic curve 
$\Lambda(\cdot)$ of arbitrary rank. Let us briefly describe 
what objects can be obtained in this way. First one can 
construct a kind of a canonical parallel transform along 
the curve $\Lambda(\cdot)$ instead of the canonical basis 
$\bigl(E_1(t),\ldots, E_m(t)\bigr)$ from Theorem 
\ref{maintheor} in the case of rank 1 curves. Namely, it 
turns out that any subspace $\Lambda(t)$ admits the 
canonical splitting 
$\Lambda(t)=\Lambda_1(t)\oplus\ldots\oplus\Lambda_s(t)$ 
such that on each subspace $\Lambda_i(t)$ the canonical 
Euclidean structure is defined $\bigl(\bigr.$for rank 1 
curves $s=m$, $\Lambda_i(t)={\rm span} 
\bigr(E_i(t)\bigr)\bigr)$; then for any $t_0$ an $t_1$ 
there exists the canonical linear mapping 
$P_{t_0,t_1}:\Lambda(t_0)\mapsto\Lambda(t_1)$ such that 
$P_{t_0,t_1}\bigl(\Lambda_i(t_0)\bigr)=\Lambda_i(t_1)$ and 
$P_{t_0,t_1}$ sends the Euclidean structure of 
$\Lambda_i(t_0)$ to the Euclidean structure of 
$\Lambda_i(t_1)$ for all $1\leq i\leq s$. Moreover, 
$P_{t_1,t_2}\circ P_{t_0,t_1}=P_{t_0,t_2}$ and 
$P_{t,t}={\rm Id}$. Secondly one can define the canonical 
complement of $\Lambda(t)$ to $W$, i.e., the subspace 
$\Lambda^{\rm comp}(t)\in L(W)$ such that 
$W=\Lambda(t)\oplus\Lambda^{\rm comp}(t)$. The main idea, 
lying in all these constructions is that if we choose some 
orthonormal basis $\bigl(e_{1i}(t),\ldots, e_{1 
m_i}(t)\bigr)$ on each subspace $\Lambda_i(t)$ (w.r.t. the 
canonical Euclidean structure on it), where $m_i=\dim 
\Lambda_i(t)$, and afterwards we take the basis on 
$\Lambda^{\rm comp}(t)$ dual (w.r.t. the symplectic form 
$\sigma$) to the basis 
$\bigl(\{e_{j1}(t)\}_{j=1}^{m_1},\bigr.$ $\ldots,$ 
$\bigl.\{e_{js}(t)\}_{j=1}^{m_s}\bigr)$ of $\Lambda(t)$, 
then the structural equation of the obtained moving Darboux 
frame in $W$ (called {\it almost canonical moving frame}) 
has to be of the simplest possible form (with the maximal 
possible trivial blocks in the matrix, corresponding to 
this structural equation). All nontrivial blocks in this 
matrix correspond to some invariant operators associated 
with our curve, which constitute a kind of the complete 
system of symplectic invariants of the curve. Note that in 
general the subspace $\Lambda^{\rm comp}(t)$, obtained in 
this way, is different from the derivative subspace 
$\Lambda^0(t)$, the construction of which is described in 
Appendix: they coincide only for regular curves. 

Finally let us describe a method for construction of 
invariants for curves in the Grassmannian $G_n(V)$ of 
$n$-dimensional subspace of $2n$-dimensional linear space 
$V$ w.r.t. the action of General Linear group $GL(V)$. It 
turns out that this problem can be reduced to the previous 
problem for the curves in Lagrange Grassmannian by an 
appropriate symplectification. Indeed, the $4n$-dimensional 
linear space $V\times V^*$ can be provided with the natural 
symplectic structure $\sigma\bigl((x_i, y_1), (x_2, 
y_2)\bigr)=y_2(x_1)-y_1(x_2)$, where $x_1, x_2\in V$ and $ 
y_1, y_2\in V^*$. 
 To any curve $\Lambda(\cdot)$ in $G_n(V)$ one can 
assign canonically the curve in Lagrange Grassmannian 
$L(V\times V^*)$. For this let $$\Lambda^{(*)}(t)={\rm 
span}\{p\in V^*: p(v)=0\,\,\, \forall v\in \Lambda(t)\}.$$ 
Then the curve $\Lambda(\cdot)\times \Lambda^{(*)}(t)$ is 
the curve in $L(V\times V^*)$. Moreover, if the curve 
$\Lambda(\cdot)$ is ample in $G_n(V)$ \footnote {The notion 
of ample curve is defined in $G_n(V)$ in the same way as in 
Lagrange Grassmannian.}, then the curve 
$\Lambda(\cdot)\times \Lambda^{(*)}(t)$ is ample. Any 
symplectic invariant of it is the invariant of the original 
curve. Besides, in this way one can construct the (almost) 
canonical moving frame also for the curves in $G_n(V)$ (in 
space $V\times V^*$). Of course, the curves obtained by the 
described symplectification are special curves in 
$L(V\times V^*)$, so some invariants from the structural 
equation of its (almost) canonical moving frame depend 
somehow one on another. 

If we start with a curve $\Lambda(\cdot)$ in the 
Grassmannian $G_k(V)$ of $k$-dimensional subspaces in $V$ 
($\dim V=2n$), where $k\neq n$, then $\Lambda(\cdot)\times 
\Lambda^{(*)}(t)$ is also the curve in $L(V\times V^*)$, 
but it is never ample. So, we cannot apply directly the 
procedure of symplectification, described above. But in 
many cases one can build from the curves in $G_k(V)$ the 
curves in $G_n(V)$ in a canonical way, combining operations 
of extension and contraction, defined by relations 
\eqref{primeicontr} and \eqref{primeiext} below, and then 
use the symplectification.

{\sl Acknowledgments:} I would like to thank professors 
Andrei Agrachev and Boris Doubrov for stimulating 
discussions.
 
\section {Proof of Theorem \ref{maintheor}}
\setcounter{equation}{0} \indent 

First let us introduce some notations. Set ${\cal 
D}^{(0)}\Lambda(t)={\cal D}_{(0)}\Lambda(t)=\Lambda(t)$ 
 and define 
inductively the following subspaces ${\cal 
D}_{(i)}\Lambda(\tau)$ and ${\cal D}^{(i)}\Lambda(\tau)$ 
for any $i\in\mathbb {N}$: 
\begin{equation}
\label{primeicontr} {\cal 
D}_{(i)}\Lambda(\tau)\stackrel{def} {=}
\left\{v\in W: 
\begin{array}{l}\exists\,\, {\rm a}\,\,{\rm curve}\,\, 
 l(\cdot)\,\,{\rm such}\,\,{\rm that}\,\,
  l(t)\in {\cal D}_{(i-1)}\Lambda(t)\,\,\forall t,\\ l(\tau)=v\,\,
  {\rm and} \,\,l^\prime(\tau)\in {\cal 
D}_{(i-1)}\Lambda(\tau)\end{array}\right \} 
\end{equation} 
\begin{equation}
\label{primeiext} {\cal D}^{(i)}\Lambda(\tau)\stackrel{def} 
{=}{\cal D}^{(i-1)}\Lambda (\tau)+ \left\{v\in W: 
\begin{array}{l}\exists\,\, {\rm a}\,\,{\rm curve}\,\, 
 l(\cdot)\,\,{\rm in}\,\,  W\,\,{\rm such}\,\,{\rm that}\,\,\\
  l(t)\in {\cal D}^{(i-1)}\Lambda(t)\,\,\forall t 
\,\,{\rm and}\,\, v= l^\prime(\tau)\end{array}\right \} 
\end{equation} 
The subspaces ${\cal D}_{(i)}\Lambda(\tau)$ and ${\cal 
D}^{(i)}\Lambda(\tau)$ are called respectively {\it the 
$i$th contraction} and {the $i$th extension} of the curve 
$\Lambda(\cdot)$ at the point $\tau$. In particular, 
directly from the definitions we have 
\begin{eqnarray}
&~&\label{ker} {\cal D}_{(1)}\Lambda(\tau)={\rm 
ker}\dot\Lambda(\tau),\\ &~& \label{rank} {\rm 
rank}\,\dot\Lambda(\tau)=\dim\, {\cal 
D}^{(1)}\Lambda(\tau)-\dim\, \Lambda(\tau). \end{eqnarray} 
Moreover, if for a given subspace $L\subset W$ we denote by 
$L^\angle$ its skew-symmetric complement, i.e. 
$L^\angle=\{v\in W: \sigma (v,l)=0\,\, \forall l\in L\}$, 
then directly from the definitions it is not hard to show 
that the subspaces ${\cal D}_{(i)}\Lambda(\tau)$ and ${\cal 
D}^{(i)}\Lambda(\tau)$ are related in the following way: 
\begin{equation} \label{subsup} {\cal D}_{(i)}\Lambda(\tau)=\Bigl({\cal 
D}^{(i)}\Lambda(\tau)\Bigr)^\angle. \end{equation} Also, 
from the definition the curve $\Lambda(\cdot)$ is ample at 
$\tau$ if and only if there exists $p\in \mathbb N$ such 
that \begin{equation} \label{amplpow} {\cal 
D}^{(p)}\Lambda(\tau)=W \,\,{\rm or}\,,\,\,{\rm 
equivalently},\,\, {\cal D}_{(p)}\Lambda(\tau)=0. 
\end{equation} 
Besides, if we suppose that the rank of $\dot\Lambda(t)$ is 
constant and equal to $r$ for any $t$, then easily 
\begin{equation}
\label{diffrank}
\begin{array}{l}
\dim\,{\cal D}^{(i)}\Lambda(t) -\dim\,{\cal 
D}^{(i-1)}\Lambda(t)\leq r,\\ \dim\,{\cal 
D}_{(i-1)}\Lambda(t) -\dim\,{\cal D}_{(i)}\Lambda(t)\leq r, 
\end{array},\quad i\in {\mathbb N}
\end{equation}

Now suppose that the curve $\Lambda(\cdot)$ has constant 
rank $1$ on some segment $I$, i.e., ${\rm 
rank}\dot\Lambda(t)=1$ for any $t\in I$. Our goal is to 
give more convenient characterization of the property of 
the rank 1 curve to be of the constant finite weight. This 
characterization is given in the following two lemmas, 
which was actually proved in \cite{jac1}. Here we 
reformulate them in our new notations. We also give a proof 
of the first of them, because it is short, while for the 
proof of the second one we refer to the corresponding 
statements from \cite{jac1}. 
\begin{lemma} {\rm (compare with Proposition 3 in 
\cite{jac1})} \label{constrank} Assume that ${\rm 
dim}\,\,W=2m$. If an ample curve $\Lambda:I\mapsto L(W)$ 
has rank 1 in the segment $I$, then out of some discrete 
subset ${\mathcal C}\in I$, one has 
\begin{equation}
\label{dimreg1} {\rm dim} {\cal 
D}^{(i)}\Lambda(t)=m+i,\quad  1\leq i\leq m,\, 
\end{equation}
or, equivalently,
\begin{equation}
\label{dimreg2} {\rm dim} {\cal 
D}_{(i)}\Lambda(t)=m-i,\quad 1\leq i\leq m,\, 
\end{equation} 
\end{lemma}

{\bf Proof.} Actually we have to prove that the set 
${\mathcal C}$ of points, where the condition 
(\ref{dimreg1}) fails, has no accumulation point. 
Otherwise, if $\bar t$ is an accumulation point of 
${\mathcal C}$, then immediately from (\ref{diffrank}) it 
follows that there are a natural number $i_0$, $ 1\leq 
i_0\leq m-1$, and a sequence of points 
$\{t_i\}_{k=1}^\infty$, converging to $\bar t$, such that 
${\cal D}^{(i_0 )}\Lambda(t_k)={\cal D}^{( 
i_0+1)}\Lambda(t_k)$ for all $k\in \mathbb N$. It implies 
that $ {\cal D}^{(j)}\Lambda(\bar t)={\cal 
D}^{(i_0)}\Lambda(\bar t)$ for any $j\geq i_0$. But by 
(\ref{diffrank}) again $\dim\, {\cal D}^{(i_0)}\Lambda(\bar 
t)\leq m+ i_0<2m$. Hence ${\cal D}^{(j)}\Lambda(\bar t)\neq 
W$ for any $j\in\mathbb{N}$. So, by (\ref{amplpow}) the 
curve $\Lambda(\cdot)$ is not ample at $\bar t$, which 
contradicts our assumptions. $\Box$ 
\begin{lemma}{\rm (see Corollary 1 and item 1 of Corollary 2 in 
\cite{jac1})} \label{cwlem} An ample curve 
$\Lambda:I\mapsto L(W)$ has the constant finite weight in 
the segment $I$ if and only if the relations 
(\ref{dimreg1}) or, equivalently, (\ref{dimreg2}) hold at 
any point of $I$. 
In this case the weight is equal to $m^2$. 
\end{lemma}
\begin{remark}
\label{flag} {\rm Actually, from the last two lemmas it 
follows that with any rank 1 curve of the constant weight 
one can associate the following curve of complete flags in 
$W$:} $$t\mapsto\Bigl({\cal 
D}_{(m-1)}\Lambda(t)\subset\ldots\subset {\cal 
D}_{(1)}\Lambda(t)\subset\Lambda(t)\subset{\cal 
D}^{(1)}\Lambda(t)\subset\ldots\subset {\cal 
D}^{(m-1)}\Lambda(t)\Bigr).\,\,\Box$$
\end{remark}

 Now let us start to prove Theorem 
\ref{maintheor}. From now one $\Lambda(\cdot)$ is a rank 1 
curve of the constant weight in $L(W)$. Without loss of 
generality it can be assumed that the velocities $\dot 
\Lambda(t)$ are nonnegative definite quadratic forms. In 
this case the curve $\Lambda(\cdot)$ is called {\it 
monotone increasing}. By the previous lemma $\dim\, {\cal 
D}_{(m-1)}\Lambda(t)=1$. For any $t$ choose a vector 
$\epsilon (t)$ such that 
\begin{equation}
\label{candir} {\cal D}_{(m-1)}\Lambda(t)={\rm 
span}\,\Bigl(\epsilon(t)\Bigr)
\end{equation}
 and the curve $t\mapsto
\epsilon (t)$ is smooth. From (\ref{primeicontr}) it 
follows easily that a smooth curve $\epsilon (\cdot)$ 
satisfies (\ref{candir}) for any $t$ if and only if the 
following relations hold 
\begin{equation} 
\label{candir1}\left\{ \begin{aligned}~&\Lambda(t)= {\rm 
span}\bigl(\epsilon(t), 
\epsilon^\prime(t),\ldots,\epsilon^{(m-1)}(t)\bigr),\\ 
~&\epsilon^{(m)}(t)\not\in\Lambda(t) 
\end{aligned}\right.\end{equation} The following lemma gives the 
canonical normalization of the vector function 
$\epsilon(t)$: 

\begin{lemma}
\label{normlem} 
There exists the unique, up to the reflection $v\mapsto 
-v$, smooth curve $\epsilon(\cdot)$ of vectors in $W$, 
satisfying (\ref{candir}), such that 
$\sigma\bigl(\epsilon^{(m)}(t), 
\epsilon^{(m-1)}(t)\bigr)=1$. 
\end{lemma}

{\bf Proof.} Let $\epsilon(\cdot)$ and 
$\tilde\epsilon(\cdot)$ be two smooth curves of vectors in 
$W$, satisfying (\ref{candir}). Then there exists a smooth 
scalar function $\alpha(t)$ such that $\tilde 
\epsilon(t)=\alpha(t)\epsilon(t)$. The last equation 
implies that 
\begin{equation}
\label{derimod}
 \tilde \epsilon^{(i)}(t)\equiv 
\alpha(t)\epsilon^{(i)}(t)\quad {\rm mod}\Bigl({\rm 
span}\bigl(\epsilon(t),\ldots\epsilon^{(i-1)}(t)\bigr)\Bigr).
\end{equation}
Note that from the first relation of (\ref{candir1}) it 
follows that \begin{equation} \label{mi0} 
\sigma(\epsilon^{(m)}(t),\epsilon^{(i)}(t))=0\quad 0\leq 
i\leq m-2. 
\end{equation}
Indeed, from the fact that all subspaces $\Lambda(t)$ are 
Lagrangian and the first relation of (\ref{candir1}) it 
follows that
\begin{equation}
\label{difss} 
\sigma(\epsilon^{(m-1)}(t),\epsilon^{(i)}(t))\equiv 0,\quad 
\sigma(\epsilon^{(m-1)}(t),\epsilon^{(i+1)}(t))\equiv 
0\quad 0\leq i\leq m-2. 
\end{equation}
Differentiating the first identity of (\ref{difss}) and 
using the second one, we obtain (\ref{mi0}). Therefore 
(\ref{derimod}) yields that 
\begin{equation}\label{srav} 
 \sigma\bigl(\tilde\epsilon^{(m)}(t), \tilde\epsilon^{(m-1)}(t)\bigr)=
 \alpha^2(t) \sigma\bigl(\epsilon^{(m)}(t), 
 \epsilon^{(m-1)}(t)\bigr).
 \end{equation}
Since by assumption $\Lambda(\cdot)$ is monotone 
increasing, one has that $$\sigma\bigl(\epsilon^{(m)}(t),
\epsilon^{(m-1)}(t)\bigr)=\dot\Lambda(t) 
\bigl(\epsilon^{(m-1)}(t)\bigr)\geq 0$$ (we use the 
identification of $\dot\Lambda(t)$ with the quadratic form, 
see (\ref{quad}); here $\dot\Lambda(t)(v)$ is the value of 
the quadratic form $\dot\Lambda(t)$ at a vector $v$). On 
the other hand, $\sigma\bigl(\epsilon^{(m)}(t), 
 \epsilon^{(m-1)}(t)\bigr)\neq 0$: Assuming the converse 
 and taking into account (\ref{mi0}), we obtain that 
 $\epsilon^{(m)}(t)\in \Lambda(t)$, 
 which contradicts  the second relation in 
 (\ref{candir1}). So, $\sigma\bigl(\epsilon^{(m)}(t), 
 \epsilon^{(m-1)}(t)\bigr)>0$. Setting $\alpha(t)=\pm\Bigl(\sigma\bigl(\epsilon^{(m)}(t), 
 \epsilon^{(m-1)}(t)\bigr)\Bigr)^{-1/2}$, we obtain from 
 (\ref{srav}) that 
 $ \sigma\bigl(\tilde\epsilon^{(m)}(t), \tilde\epsilon^{(m-1)}(t)\bigr)=1.$
 It remains only to notice that by our constructions  
 $\tilde\epsilon(t)$, satisfying the last relation, is defined up to 
 the sign. $\Box$
 
 Now suppose that $\epsilon(t)$ is as in the previous lemma. 
 We set
 \begin{equation}
 \label{Ecan}
 E_i(t)=\epsilon^{(i-1)}(t), \quad 1\leq i\leq m;\quad 
 F_m(t)=\epsilon^{(m)}(t)
 \end{equation}
 
By the first relation of (\ref{candir1}) the tuple 
$\bigl(E_1(t),\ldots, E_m(t)\bigr)$ satisfies the first 
condition of Theorem \ref{maintheor}, while together with 
$F_m(t)$ it satisfies the first two equations of 
(\ref{structeq}). Besides, by (\ref{mi0}) the choice of the
vectors, defined in (\ref{Ecan}), does not contradict the 
relations (\ref{darboux}). 

To finish the proof of the theorem it remains to complete 
the tuple $\bigl((E_1(t),\ldots, E_m(t),F_m(t)\bigr)$ to 
the moving Darboux frame in $W$, which satisfies 
the last two equations of (\ref{structeq}), and to show 
that such complement is unique (the freedom in the sign, 
mentioned in the last sentence of Theorem \ref{maintheor} 
will follow then from the freedom up to the sign in the 
choice of $\epsilon(t)$ in Lemma \ref{normlem}). 
   
For this we analyze the structural equations of all 
possible moving Darboux's frames, and choose among them  
one, which has the maximal possible number of zero entries 
in the matrix of its structural equation. First take some 
tuple $\{\overline F_i(t)\}_{i=1}^{m-1}$ such that 
\begin{equation}
\label{frame1} \bigl(E_1(t),\ldots, E_m(t), \overline 
F_1(t), \ldots,\overline F_{m-1}(t), F_m(t)\bigr) 
\end{equation}
is a moving Darboux's frame in $W$. Then from the 
definition of Darboux's basis (see (\ref{darboux})), 
and the first two equations of (\ref{structeq}) it follows 
that there exist functions $\bar \xi_{ij}(t)$, $1\leq 
i,j\leq m-1$ such that 
\begin{equation}
\label{structeq1} \overline 
{F}_i^{\,\prime}(t)=\sum_{j=1}\bar \xi_{ij}(t)E_j(t)- 
(1-\delta_{1i})\overline F_{i-1}(t),\quad 
\bar\xi_{ij}(t)=\bar\xi_{ji}(t),\quad 1\leq i,j\leq m-1, 
\end{equation}
where $\delta_{kl}$ is the Kronecker symbol. Further, by 
the definition of Darboux's basis for the given tuple 
$\{\widehat F_i(t)\}_{i=1}^{m-1}$ of curves of vector in 
$W$ the frame 
\begin{equation}
\label{frame2} \bigl(E_1(t),\ldots, E_m(t), \widehat 
F_1(t), \ldots,\widehat F_{m-1}(t), F_m(t)\bigr) 
\end{equation}
is a moving Darboux's frame in $W$ if and only if there 
exist functions $b_{ij}(t)$, $1\leq i,j\leq m-1$ such that 
\begin{equation}
\label{transDb} \left\{\begin{aligned}~&\widehat 
F_i(t)=\overline F_i(t)+ \sum_{j=1}^{m-1}b_{ij}(t)E_j(t), 
\quad 1\leq i\leq m-1\\~ & b_{ij}(t)=b_{ji}(t) , \quad 
1\leq i,j\leq m-1. 
\end{aligned}
\right. 
\end{equation} 
Besides, similarly to (\ref{structeq1}), for the tuple 
$\{\widehat F_i(t)\}_{i=1}^{m-1}$ there exist functions 
$\hat\xi_{ij}(t)$, $1\leq i,j\leq m-1$ such that 
\begin{equation}
\label{structeq2} \widehat F_i^{\,\prime}(t)= 
\sum_{j=1}\hat\xi_{ij}(t)E_j(t)-(1-\delta_{1i})\widehat 
F_{i-1}(t),\quad \hat\xi_{ij}(t)=\hat\xi_{ji}(t),\quad  
1\leq i,j\leq m-1. 
\end{equation}
Now we are ready to find the transformation rule from the 
coefficients $\bar\xi_{ij}$ of the structural equation for 
the original frame (\ref{frame1}) to the coefficients 
$\hat\xi_{ij}$ of the structural equation for the frame 
(\ref{frame2}): Set 
\begin{equation}
\label{bm} b_{im}(t)=b_{mi}(t)=0,\quad 1\leq i\leq m. 
\end{equation}
Then, substituting (\ref{transDb}) into (\ref{structeq2}), 
and using (\ref{structeq1}) one can easily obtain 
\begin{equation}
\label{transxi}
\hat\xi_{ij}(t)=\bar\xi_{ij}(t)+b_{ij}^\prime(t)+
(1-\delta_{1i}) b_{i-1, j}(t)+(1-\delta_{1j}) b_{i, 
j-1}(t). 
\end{equation}
From transformation rule (\ref{transxi}) it follows 
immediately that Theorem \ref{maintheor} will follow from 
the following 
\begin{lemma}
\label{diaglem} For the given smooth curve of $m\times m$ 
symmetric matrices 
$\overline\Omega(t)=\bigl(\bar\xi_{ij}(t)\bigr)_{ij=1}^m$ 
there exists the unique smooth curve of symmetric $m\times 
m$ matrices $B(t)=\bigl(b_{ij}(t)\bigr)_{ij=1}^m$, 
satisfying (\ref{bm}), such that 
\begin{equation}
\label{transxi0} 
\bar\xi_{ij}(t)+b_{ij}^\prime(t)+ (1-\delta_{1i}) b_{i-1, 
j}(t)+(1-\delta_{1j}) b_{i, j-1}(t)=0,\quad i\neq j,
\end{equation} 
where $\delta_{kl}$ is the Kronecker symbol. 
\end{lemma}
{\bf Proof.} For $m=1$ there is nothing to prove. Suppose 
that $m>1$. Note that by the symmetricity of equations 
(\ref{transxi0}) w.r.t. permutation $(ij)\mapsto(ji)$ it is 
enough to prove existence and uniqueness of $b_{ij}(t)$ 
with $i\geq j$. We will "fill" step by step the lower 
triangle (including the diagonal) of the matrix $B(t)$ 
starting from the $(m-1)$th row (the last row is given by 
(\ref{bm})). Taking into account (\ref{bm}), from equation 
(\ref{transxi0}) for $i=m$ it follows that 
$$b_{m-1, j}(t)=-\bar\xi_{mj}(t),\quad 1\leq j\leq m-1.$$ 
In this way and using symmetricity we have filled $(m-1)$th 
row of $B(t)$. 

Now suppose by induction that for some $2<i\leq m-1$ we 
have filled $\bar i$th rows of the matrix $B(t)$ for all 
$\bar i\geq i$ (note that if $i=2$, we are already done 
from symmetricity). We would like to determine the 
$(i-1)$th row. 
From equation (\ref{transxi0}) for $j=1$ it follows that 
\begin{equation}
\label{m2row} b_{i-1,1}(t)=-\bar\xi_{i, 1}(t)-b_{i, 
1}^\prime(t). \end{equation} Since the righthand side of 
(\ref{m2row}) is determined by the induction hypothesis, 
$b_{i-1, 1}(t)$ is determined. Other elements $b_{i-1, 
j}(t)$ with $2\leq j\leq i$ are determined from the 
following recursive formula $$ b_{i-1, j}(t)=-\bar\xi_{i 
j}(t)-b_{i j}^\prime(t)-b_{i-1, j-1}(t),$$ which follows 
from (\ref{transxi0}). In this way and using symmetricity 
we have filled $(i-1)$th row of $B(t)$. The proof by 
induction is completed.$\Box$ 

Let $B(t)=\bigl(b_{ij}(t)\bigr)_{ij=1}^m$ be as in the 
previous lemma. Then, setting $F_i(t)=\widehat F_i(t)$, 
$1\leq i\leq m-1$, where $\widehat F_i(t)$ is defined by 
the first equation of (\ref{transDb}), we obtain the moving 
Darboux frame, required in Theorem \ref{maintheor}. Note 
that from (\ref{transxi}) with $j=i$ it follows that the 
functions $\lambda_{i}$ from (\ref{structeq}) satisfy 
$$\lambda_{m-i+1}(t)=\bar\xi_{ii}(t)+b_{ii}^\prime(t)+ 
2(1-\delta_{1i}) b_{i, i-1}(t), 
\quad 1\leq i\leq m.$$ The proof of Theorem 
\ref{maintheor} is completed. $\Box$ 

\section{Proof of Proposition \ref{equivpr}}
\setcounter{equation}{0}\indent 

Throughout this section we will use the following 
notations: For a given tuple $\{\psi_i(t)\}_{i=1}^N$ of 
smooth functions we denote by ${\rm 
Pol}\bigl(\{\psi_i(t)\}_{i=1}^N\bigr)$ any function, which 
can be expressed as a polynomial (over $\mathbb R$) without 
free term w.r.t. the functions $\psi_i(t)$, $1\leq i\leq 
N$, and their derivatives. Also, we denote by ${\rm 
Lin}\bigl(\{\psi_i(t)\}_{i=1}^N\bigr)$ any function, which 
can be expressed as a linear combination (over $\mathbb R$) 
of the functions $\psi_i(t)$, $1\leq i\leq N$, and their 
derivatives. 

One can try to prove Proposition \ref{equivpr}, expressing 
$\beta_{2i}(t)$, $0\leq i\leq m-1$, by $\lambda_j(t)$ with 
the help of the structural equation (\ref{structeq}) and 
formula (\ref{gdet}): Let $W\cong\mathbb R^m\times\mathbb 
R^m$,$\Lambda(t)=\{(x, S_t x):x\in \mathbb R^m\}$, 
$E_i(t)=\bigl(\pi_{1i}(t),\ldots,\pi_{2m,i}(t)\bigr)$. 
Denote also by $\Pi_1(t)$ and $\Pi_2(t)$ the following 
$m\times m$-matrices: $$\Pi_1(t)=(\pi_{ji}(t))_{1\leq j\leq 
m, 1\leq i\leq m},\quad \Pi_2(t)=(\pi_{ji}(t))_{m+1\leq 
j\leq 2m, 1\leq i\leq m}$$ Then 
$S_t=\Pi_2(t)\Pi_1(t)^{-1}$. Using the structural equation 
(\ref{structeq}) one can find derivatives of $S_t$ of any 
order. So, using (\ref{gdet}), one can compute in general 
the Taylor formula for $g_{_{\Lambda}}$ up to the required 
order. But in this way one meets rather cumbersome 
computations even in the case $m=2$. 

Our proof of Proposition \ref{equivpr} basically consists 
of the following two steps: 

{\bf Step 1.} From (\ref{structeq}) one can express 
$E_1^{(2m)}(t)$ as a linear combination of $E_1(t),\ldots, 
E_1^{(2m-1)}(t)$: it turns out that if  
$E_1^{(2m)}(t)=\displaystyle {\sum_{k=0}^{2m-1}}\Gamma_k(t) 
E_1^{(k)}(t)$, then  
\begin{equation}
\label{coef1} \Gamma_{2i}(t)=(-1)^i\lambda_{m-i}(t)+ {\rm 
Lin}\bigl(\{\lambda_j(t)\}_{j=1}^i\bigr),\quad 0\leq i\leq 
m-1 
\end{equation} 
Indeed, using the first, second, and forth relations of 
(\ref{structeq}), it is easy to show by induction that for 
all $1\leq s\leq m-1$ 
\begin{equation}
\label{ind1} \begin{split} &F_{m-s}(t)=(-1)^s 
E_1^{(m+s)}(t)+\sum_{k=1}^{s-1}\Bigl((-1)^{s-k}\lambda_k(t)+{\rm 
Lin}\bigl(\{\lambda_l(t)\}_{l=1}^{k-1}\bigr)\Bigr)E_1^{(m+s-2k)}(t)+\\ 
&\sum_{k=1}^{s-1}\Bigl((-1)^{s-k}(s-k)\lambda_k^\prime(t)+{\rm 
Lin}\bigl(\{\lambda_l(t)\}_{l=1}^{k-1}\bigr)\Bigl) 
E_1^{(m+s-2k-1)}(t)+\lambda_s(t) E_1^{(m-s)}(t). 
\end{split} \end{equation} 
Then substituting (\ref{ind1}) for $s=m-1$ in the third 
relation of (\ref{structeq}) one gets (\ref{coef1}).

{\bf Step 2.} First let us give a sketch of what we are 
going to do. Let 
\begin{equation}
\label{oldbasis} \bigl(e_1(t),\ldots, e_m(t),f_1(t),\ldots 
f_m(t)\bigr)
\end{equation}
 be the canonical basis of the curve 
$\Lambda(\cdot)$, constructed in \cite{jac1} 
(for the sake of completeness we will describe this 
construction in Appendix). Collecting some information 
about its structural equation from \cite{jac1} and 
\cite{jac2}, we will express $e_1^{(2m)}(t)$ as a linear 
combination of $e_1(t),\ldots e_1^{(2m-1)}(t)$. Namely, if 
$e_1^{(2m)}(t)=\displaystyle{\sum_{k=0}^{2m-1}}\gamma_k(t) 
e_1^{(k)}(t)$, then 
\begin{equation}
\label{coef2} \gamma_{2i}(t)=\overline C_i 
\beta_{2(m-1-i)}(t)+ {\rm 
Pol}\bigl(\{\beta_{2j}(t)\}_{j=0}^{i-1}\bigr),\quad 0\leq 
i\leq m-1. 
\end{equation} 
 On the other hand, from the structural equation again 
 it will follow that 
 \begin{equation}
 \label{newold}
 e_1(t)={\rm const}\, E_1(t),
 \end{equation}
which implies that $\gamma_k(t)={\rm const}\, \Gamma_k(t)$. 
Comparing (\ref{coef1}) with (\ref{coef2}) in view of the 
last relation, one gets (\ref{equivbl}), which together 
with Remark \ref{mpoint} implies Proposition \ref{equivpr}. 

Now we start to prove formulas (\ref{coef2}) and 
(\ref{newold}). 
All information from \cite{jac1} and \cite{jac2} that we 
need about the frame (\ref{oldbasis}) can be summarized in 
the following 
\begin{lemma}
\label{structlem} The frame (\ref{oldbasis}), constructed 
in \cite{jac1}, satisfies the following equation 
\begin{equation}
\label{oldstruct}
 \left\{\begin{aligned} 
~&e_i'(t)=\sum_{j=1}^m\alpha_{ij}(t)e_j(t)+m^2\delta_{mi}f_m(t)\\
~&f_i'(t)=\sum_{j=1}\sigma_{ij}(t)e_j(t)-\sum_{j=1}^m\alpha_{ji}(t)f_j(t),
\end{aligned}\right. 
\end{equation}
where $\sigma_{ij}(t)=\sigma_{ji}(t)$; $\delta_{mi}$ is the 
Kronecker symbol;
\begin{eqnarray}
&~&\label{i1} \alpha_{ij}(t)\equiv 0 \,\, {\rm for}\,\, 
i<j-1;\quad \alpha_{i-1,i}(t)\equiv 
\frac{(i-1)(2m-i+1)}{m-i+1};\\ 
&~&\label{i3}\alpha_{ii}=0;\\ &~& 
\label{i4}\alpha_{ij}(t)=\left\{\begin{array}{ll} 
\nu_{ij}\beta_{i-j-1}(t)+{\rm 
Pol}\Bigl(\{\beta_{2s}(t)\}_{s=0}^{^{\frac{i-j-3}{2}}}\Bigr),& 
i-j\,\, {\rm is}\,\, {\rm positive}\,\,{\rm odd,} 
\\ ~\\{\rm 
Pol}\Bigl(\{\beta_{2s}(t)\}_{s=0}^{^{\frac{i-j-2}{2}}}\Bigr),& 
i-j\,\, {\rm is}\,\, {\rm positive}\,\,{\rm even}; 
\end{array}\right.\\
&~&\label{i5} \sigma_{ij}(t)=\left\{\begin{array}{ll} 
c_{ij}\beta_{2m-i-j}(t)+{\rm 
Pol}\Bigl(\{\beta_{2s}(t)\}_{s=0}^{^{m-\frac{i+j+2}{2}}}\Bigr),& 
i+j\,\, {\rm is}\,\,{\rm even,} 
\\ ~\\{\rm 
Pol}\Bigl(\{\beta_{2s}(t)\}_{s=0}^{^{m-\frac{i+j+1}{2}}}\Bigr),& 
i+j\,\, {\rm is}\,\,{\rm odd}; 
\end{array}\right.
\end{eqnarray} 
$\nu_{ij}$ and $c_{ij}$ are some constants.
\end{lemma}

Relations (\ref{i1}) are exactly items 1 and 2 of Lemma 7.3 
from \cite{jac1}. Relation (\ref{i3}) for $1\leq i\leq m-1$ 
is exactly relation (1.79) in \cite{jac2}, while for $i=m$ 
it can be obtained easily from formula (7.52) of 
\cite{jac1}, taking into account identity (1.74) from 
\cite{jac2}. Further, relation (\ref{i4}) for $1\leq i\leq 
m-1$ is more specified version of item 3 of Lemma 7.3 from 
\cite{jac1}, which follows from the proof of this lemma, 
while for $i=m$ it can be obtained without difficulties 
from formula (7.52) of \cite{jac1}. Finally, relation 
(\ref{i5}) in the case of even $i+j$ is exactly Lemma 1.6 
of \cite{jac2}, while in the case of odd $i+j$ it follows 
from the proof of this lemma (see, for example, formula 
(1.67) there). 

One can reformulate Lemma \ref{structlem} in more 
convenient form if one denotes 
\begin{equation}
\label{convert} \forall\, 1\leq i\leq m:\quad 
e_{m+i}(t)\stackrel{def}{=}f_{m-i+1}(t);\quad {\mathcal 
E}(t)\stackrel{def}{=}\begin{pmatrix} e_1(t)\\ 
\vdots\\e_{2m}(t)\end{pmatrix} .\end{equation} If 
${\mathcal M}(t)$ is $2m\times 2m$ matrix, ${\mathcal 
M}(t)=\{\mu_{ij}(t)\}_{i,j=1}^{2m}$, such that 
\begin{equation}
\label{oldstruct1}
 {\mathcal E}^\prime(t)= {\mathcal M}(t) 
{\mathcal E}(t),
\end{equation}
 then from (\ref{oldstruct})-(\ref{i5}) it 
follows easily that 
\begin{equation}
\label{ni} 
\mu_{ij}(t)=\left\{\begin{array}{ll}0&j>i+1\,\,{\rm 
and}\,\,j=i\\ \chi_{ij}& j=i+1 
\\\chi_{ij}\beta_{i-j-1}(t)+{\rm 
Pol}\Bigl(\{\beta_{2s}(t)\}_{s=0}^{^{\frac{i-j-3}{2}}}\Bigr),& 
i-j\,\, {\rm is}\,\, {\rm positive}\,\,{\rm odd,} 
\\ ~\\{\rm 
Pol}\Bigl(\{\beta_{2s}(t)\}_{s=0}^{^{\frac{i-j-2}{2}}}\Bigr),& 
i-j\,\, {\rm is}\,\, {\rm positive}\,\,{\rm even},
\end{array}\right.
\end{equation}
where $\chi_{ij}$ are constants. Set also $\chi_{2m, 
2m+1}=1$. Then combining (\ref{oldstruct1}) with 
(\ref{ni}), one can obtain without difficulties by 
induction that 
\begin{equation} \label{induct}(1- \delta_{i, 2m}) 
e_{i+1}(t)=\Bigl(\displaystyle{\prod_{k=1}^i\chi_{k,k+1}}\Bigr)^{-1}
\Bigl(\sum_{j=0}^{i-2} \kappa_{ij}(t)e_1^{(j)}(t) 
+e_1^{(i)}(t)\Bigr), \quad 0\leq i\leq 2m-1, 
\end{equation} where $\delta_{i,2m}$ is the Kronecker 
symbol,
\begin{equation}
\label{kappa} \kappa_{ij}(t)=\left\{\begin{array}{ll} 
\rho_{ij}\beta_{i-j-2}(t)+ {\rm 
Pol}\Bigl(\{\beta_{2s}(t)\}_{s=0}^{^{\frac{i-j-4}{2}}}\Bigr)&i-j\,\, 
{\rm is}\,\, {\rm positive}\,\,{\rm even}\\ 
{\rm 
Pol}\Bigl(\{\beta_{2s}(t)\}_{s=0}^{^{\frac{i-j-3}{2}}}\Bigr)&i-j\,\, 
{\rm is}\,\, {\rm positive}\,\,{\rm odd} 
\end{array}\right.,
\end{equation}
\begin{equation}
\label{compform} 
\rho_{ij}=\displaystyle{\sum_{s=1}^{j+1}}\left(\chi_{i-j-1+s,s} 
\displaystyle{\prod_{l=s}^{i-j-2+s}}\chi_ 
{l,l+1}\right),\quad i-j\,\, {\rm is}\,\, {\rm 
positive}\,\,{\rm even}. 
\end{equation}

Relation (\ref{induct}), used for $i=2m$, implies 
\eqref{coef2}. Actually, $\overline C_i$ in (\ref{coef2}) 
can be taken as \begin{equation} \label{overallc} \overline 
C_i=-\rho_{2m, i} 
\Bigl(\displaystyle{\prod_{k=1}^i\chi_{k,k+1}}\Bigr)^{-1}. 
\end{equation} 
Further, relation (\ref{induct}), used for $i\leq m$, 
implies that $e_1(t)$, taken as $\epsilon(t)$, satisfies 
\eqref{candir1}. Therefore $e_1(t)=\alpha(t) E_1(t)$. 
Moreover, $$\sigma\bigl(e_1^{(m)},e_1^{(m-1)})= 
\chi_{m,m+1}\displaystyle{\prod_{l=1}^{m-1}}\chi_{l,l+1}^2,$$ 
which implies that $\alpha(t)$ is constant (equal to 
$(\chi_{m,m+1})^{1/2}\displaystyle{\prod_{l=1}^{m-1}} 
\chi_{l,l+1}$ ). This proves formula (\ref{newold}), which 
completes the proof of Proposition \ref{equivpr}. $\Box$ 

\begin{remark}
\label{diffproof} {\rm Actually, the existence part of 
Theorem \ref{uniex1} will follow from (\ref{coef2}) only ( 
without using the modified principal curvatures 
$\lambda_j(t)$ and Proposition \ref{equivpr}), if one shows 
that \begin{equation} \label{Cn0} 
 \overline C_i\neq 0,\quad 0\leq i\leq m-1.
 \end{equation}
 But even after finding explicit expressions for $\chi_{ij}$ 
 we did not succeed  to verify (\ref{Cn0}) , using 
 formula \eqref{overallc}. Instead, the presence of the 
 complete system of invariants $\{\lambda_i(t)\}_{i=1}^m$ and 
 identity (\ref{equivbl}) imply (\ref{Cn0}) 
 automatically, as was mentioned already in Remark \ref{mpoint} above. 
 } 
\end{remark}
\section{Appendix}
\setcounter{equation}{0} \indent 

In this Appendix we briefly describe the construction of 
the canonical moving Darboux frame for a rank 1 curve of 
the constant weight, introduced in \cite{jac1} and used in 
the previous section. 
The construction is based on the fact that the set 
$\Lambda^\pitchfork$ of all Lagrangian subspaces 
transversal to a subspace $\Lambda\in L(W)$ can be 
naturally endowed with the structure of an affine space 
over the linear space ${\rm Quad}\, (\Lambda^*)$ of 
quadratic forms on $\Lambda^*$ or, equivalently, over the 
linear space ${\rm Symm}\,(\Lambda^*)$ of self-adjoint 
linear mappings from $\Lambda^*$ to $\Lambda$. Indeed, as 
in Introduction, for a given $\Gamma \in 
\Lambda^\pitchfork$ denote by ${\mathcal 
B}_{\Lambda,\Gamma}$ the following linear mapping from 
$\Gamma$ to $\Lambda^*$ : $v\mapsto \sigma(v, \cdot), \quad 
v\in \Lambda_1$, $v\in\Gamma$. Then the operation of 
subtraction on $\Lambda^\pitchfork$ with values in ${\rm 
Quad} (\Lambda^*)$ can be defined as follows : 
\begin{equation} \label{subtr} 
(\Gamma-\Delta)(l)=\sigma\bigl(({\mathcal 
B}_{\Lambda,\Gamma})^{-1}l,({\mathcal 
B}_{\Lambda,\Delta})^{-1} l),\quad \Gamma,\Delta \in 
\Lambda^\pitchfork,\,\, l\in\Lambda^* .
\end{equation}
It is not difficult to show that $\Lambda^\pitchfork$ 
endowed with this operation of subtraction satisfies the 
axioms of affine space.

Consider now some curve $\Lambda(\cdot)$ in $L(W)$. Fix 
some parameter $\tau$. Note that if the curve 
$\Lambda(\cdot)$ is ample at $\tau$, then 
$\Lambda(t)\in\Lambda(\tau)^\pitchfork$ for all $t$ from a 
punctured neighborhood of $\tau$. Then we obtain the curve 
$t\mapsto\Lambda(t)\in\Lambda(\tau)^\pitchfork$ in the 
affine space $\Lambda(\tau)^\pitchfork$. Denote by 
$\Lambda_\tau(t)$ the identical embedding of $\Lambda(t)$ 
in the affine space $\Lambda(\tau)^\pitchfork$. The 
velocity $\frac{\partial}{\partial t}\Lambda_\tau(t)$ is an 
element of the underlying linear space, i.e., it is well 
defined self-adjoint mappings from $\Lambda^*$ to 
$\Lambda$. Now let $\Lambda(\cdot)$ be a rank $1$ curve in 
$L(W)$. For definiteness suppose that it is monotone 
increasing. Then $\frac{\partial}{\partial 
t}\Lambda_\tau(t)$ is a nonpositive definite rank $1$ 
self-adjoint linear mapping from $\Lambda^*$ to $\Lambda$ 
and for $t\neq \tau$ there exists a unique, up to the sign, 
vector $w(t,\tau) \in \Lambda(\tau)$ 
  such that
$\langle v,\frac{\partial}{\partial t}\Lambda_\tau(t) 
v\rangle= - \langle v,w(t,\tau)\rangle^2$ for any $v\in 
\Lambda(\tau)^*$. 
The properties of the vector function 
$t\mapsto w(t,\tau)$ 
can be summarized as follows ( see \cite{jac1}, 
section 7, Proposition 4 and Corollary 2): 

\begin{prop}
\label{constcor} If $\Lambda(\cdot)$ is a rank 1 curve of 
the constant weight in $L(W)$, then for any $\tau$ the 
function $t\mapsto w(t,\tau)$ has a pole of order $m$ at 
$t=\tau$. Moreover, if we write down the expansion of 
$t\mapsto w(t,\tau)$ in Laurent series at $t=\tau$, 
\begin{equation} \label{weq} 
w(t,\tau)=\sum_{i=1}^{m}e_i(\tau)(t-\tau)^{i-1-l}+O(1),
\end{equation}
then the vector coefficients $e_1(\tau),\ldots,e_m(\tau)$ 
constitute a basis of the subspace $\Lambda(t)$. 
\end{prop}

So, formula (\ref{weq}) defines the canonical basis 
$e_1(\tau),\ldots, e_m(\tau)$ on each subspace 
$\Lambda(\tau)$ of the rank $1$ curve $\Lambda(\cdot)$ of 
constant weight. In order to complete this basis to some 
canonical moving Darboux's frame in $W$, one can exploit 
the affine structure again: 
Fixing an 
``origin'' $\Delta$ in $\Lambda(\tau)^\pitchfork$ we obtain 
 a vector function $t\mapsto \bigl(\Lambda_\tau(t)-\Delta\bigr)$ 
 with values in ${\rm 
Quad}\,(\Lambda^*)$ (or ${\rm Symm}\,(\Lambda^*)$). 
Actually, the fact that the curve $\Lambda(\cdot)$ is ample 
at the point $\tau$ is equivalent to the fact that the 
vector function $t\mapsto 
\bigl(\Lambda_\tau(t)-\Delta\bigr)$ has the pole at 
$t=\tau$. Using only the axioms of affine space, one can 
prove easily that there exists a unique subspace 
$\Lambda^0(\tau)\in\Lambda^\pitchfork$ such that the free 
term in the expansion of the vector function $t\mapsto 
\bigl(\Lambda_\tau(t)-\Lambda^0(\tau)\bigr)$ to the Laurent 
series at $\tau$ is equal to zero. The curve 
$\tau\mapsto\Lambda^0(\tau)$ is called the {\it derivative 
curve of the ample curve $\Lambda(\cdot)$}. 
Now let $f_1(\tau),\ldots, f_m(\tau)$ be a basis of 
$\Lambda^0(\tau)$ dual to the canonical basis of 
$\Lambda(\tau)$, i.e. $\sigma (f_i(\tau), 
e_j(\tau))=\delta_{ij}$. The tuple 
$\bigl(e_1(\tau),\ldots,e_m(\tau),f_1(\tau),\ldots,f_m(\tau)\bigr)$ 
is exactly the canonical moving Darboux's frame in $W$, 
properties of which we used in section 3.


\end{document}